\documentclass[12pt,a4paper]{article}
\usepackage{amsmath,amsthm,amssymb,latexsym,epic}
\usepackage{graphicx}

\newtheorem{theorem}{Theorem}
\newtheorem{lemma}{Lemma}

\newtheorem{proposition}{Proposition}

\usepackage[all]{xy}
\usepackage{enumerate}

\renewcommand{\ker}{\mathrm{ker}}

\newcommand{\ISA}{{\mathcal {ISA}}}

\newcommand{\IS}{{\mathcal {IS}}}

\begin{document}
\markboth{G. Kudryavtseva} {The structure of automorphism groups
of semigroup inflations}

%\catchline{}{}{}{}{}
%
%%%%%%%%%%%%%%%%%%%%%%%%%%%%%%%%%%%%%%%%%%%%%%%%%%%%%%%%%%%%%%%%%%%%

\title{THE STRUCTURE OF AUTOMORPHISM GROUPS OF SEMIGROUP INFLATIONS }
\author{GANNA KUDRYAVTSEVA}
 \date{}
 \maketitle
\begin{center}{Department of Mechanics and Mathematics, \\
Kyiv Taras Shevchenko University,\\ Volodymyrs'ka, 60, Kuiv, 01033, Ukraine\\
akudr@univ.kiev.ua}
\end{center}

\begin{abstract}
It is proved that the automorphism group of a semigroup being an
inflation of its proper subsemigroup decomposes into a semidirect
product of two groups one of which is a direct sum of full
symmetric groups.
\end{abstract}

\section{Introduction}\label{s1}

In the study of a specific semigroup the description of all its
automorphisms is one of the most important questions. The
automorphism groups of many important specific semigroups are
described (see, for example, the references in ~\cite{Sz}). It
happens that for two types of examples of semigroups having rather
different nature the automorphism groups have similarities in
their structure: each of them decomposes into a semidirect product
of two groups one of which is a direct sum of the full symmetric
groups. These two types of semigroups are variants of some semigroups
of mappings (see~\cite{Symons}, \cite{KTs}) and maximal nilpotent
subsemigroups of some transformation semigroups (see~\cite{GTS},
\cite{Str}).

In the present paper we give a general setting for all
these results via consideration of the automorphism group of an
abstract semigroup being an inflation of its proper subsemigroup, constructed
as follows.  Define an equivalence $h$ on a
semigroup $S$ via $(a,b)\in h$ if and only if $ax=bx$ and $xa=xb$
for all $x\in S$, i.e. $a$ and $b$ are not distinguished by
multiplication from either side. The notation $h$ goes back
to~\cite{Symons}, where $h$ was defined and studied for the
variants of certain semigroups of mappings. Set further
$$\psi=(h\cap((S\setminus S^2)\times (S\setminus S^2)))\cup \{(a,a):a\in S^2\}.$$
Denote by $T$ any transversal of $\psi$.
Then $T$ is a subsemigroup of $S$, and $S$ is an inflation of $T$.
An automorphism $\tau$ of $T$ will be called {\em extendable}
provided that $\tau$ coincides with the restriction to $T$ of a
certain automorphism of $S$. Clearly, all extendable automorphisms
of $T$ constitute a subgroup, $H$, of the group $Aut T$ of all
automorphisms of $T$.

Our main result is the following theorem.

\begin{theorem} \label{th:main} The group ${\mathrm{Aut}} S$ is isomorphic
to a semidirect product of two groups one of which (the one which is normal) is the
direct sum of the full symmetric groups on the $\psi$-classes
and the other one is the group $H$ consisting of all extendable automorphisms
of $T$.
\end{theorem}

In all the papers \cite{Symons}, \cite{KTs}, \cite{GTS},
\cite{Str}  the semigroups under consideration are  inflations on some their proper
subsemigroups (although this is not mentioned explicitly in any of these papers),
and one of the multiples (the one which is normal) in the structure theorem for automorphism group
every time coincides with the corresponding multiple from Theorem~\ref{th:main}
i.e. with the direct sum of the full symmetric groups on $\psi$-classes.
Only the structure of the other multiple depends on the specific of
a semigroup under consideration. However, in all the papers \cite{Symons}, \cite{KTs}, \cite{GTS},
\cite{Str} the structure of this multiple is easily understood, and the main points each time
were to show that these easily understood automorphisms do not exhaust all the automorphisms
of a given semigroup, to find the normal multiple from Theorem¨\ref{th:main} and to show
that the construction of a semidirect product comes to the game.

\section{Construction}

A semigroup $S$ is called an {\em inflation} of its subsemigroup
(see \cite{clifford},  section 3.2)
$T$ provided that there is an onto map $\theta:S\to T$ such that:
\begin{itemize}
\item $\theta^2=\theta$;
\item $a\theta b\theta=ab$ for all $a,b\in S$.
\end{itemize}
In the described situation $S$ is often referred to as an {\em
inflation of} $T$ {\em with an associated
map} $\theta$ (or just {\em with a map }$\theta$).

It is immediate that if $S$
is an inflation of $T$ then $T$ is a retract of $S$ (i.e. the image under
an idempotent homomorphism) and that $S^2\subset T$.

\begin{lemma} \label{lem:mon} Suppose that $S$ is an inflation of $T$ with a map
$\theta$. Then $\ker\theta\subset h$.
\end{lemma}

\begin{proof}
Let $(a,b)\in\ker\theta$ and $s\in S$. Then we
have
\begin{equation*}
as=a\theta s\theta = b\theta s\theta =bs;\, sa=s\theta a\theta
=s\theta b\theta=sb.
\end{equation*}
It follows that $(a,b)\in h$.
\end{proof}

\begin{lemma}
The equivalence $\psi$, defined in the Introduction, is a congruence on $S$.
\end{lemma}
\begin{proof}
Obviously, $\psi$ is an equivalence relation.
Prove that $\psi$ is left and right compatible.  Let $(a,b)\in\psi$ and $a\neq b$. Then
$(a,b)\in (h\cap((S\setminus S^2)\times (S\setminus S^2)))$.
Let $c\in S$. As $(a,b)\in h$, one has $ac=bc$ and $ca=cb$ for each $c\in S$.
It follows that $(ac,bc)\in\psi$ and $(ca,cb)\in\psi$ as $\psi$ is reflexive.
\end{proof}

Set $T$ to be a transversal of $\psi$.
\begin{lemma}
$T$ is a subsemigroup of $S$, and $S$ is an inflation of $T$.
\end{lemma}

\begin{proof}
$T$ is a subsemigroup of $S$ as $T\supset S^2$. Let $\theta$ be the
map $S\to T$ which sends any element $x$ from $S$ to the unique element
of the $\psi$-class of $x$, belonging to $T$. The construction implies that $S$ is an
inflation of  $T$ with the map $\theta$.
\end{proof}

 Let
$S=\cup_{a\in T}X_a$ be the decomposition of $S$ into the union of
$\psi$-classes, where $X_a$ denotes the $\psi$-class of $a$.
Set $G_a$ to be the full symmetric group acting on $X_a$ and
$G=\oplus_{a\in T}G_a$.

We start from the following easy observation

\begin{lemma}
$\pi$ is an automorphism of $S$ for each $\pi\in G$.
\end{lemma}

\begin{proof}
It is enough to show that $(xy)\pi=x\pi y\pi$ whenever $x,y\in S$.
Suppose first that  $x, y\in S\setminus S^2$. Since $xy\in S^2$ it
follows that $\pi$ stabilizes $xy$, so that $(xy)\pi=xy$. Now, the
inclusions $(x, x\pi)\in h$ and $(y,
y\pi)\in h$ imply $x\pi y\pi =x\pi y = xy$. This
yields $xy\pi=x\pi y\pi$, and the proof is complete.
\end{proof}

The following proposition gives a characterization of extendable
automorphisms of $T$.

\begin{proposition}\label{lem:thues} An automorphism $\tau$ of $T$ is extendable if and
only if the following condition holds:
\begin{equation}\label{eq:q}
(\forall a,b\in T)\,\,\, a\tau= b \,\Rightarrow \,|X_a|=|X_b|.
\end{equation}
\end{proposition}

\begin{proof}Suppose $\tau\in{\mathrm{Aut}}T$ is extendable and $a\in T$.
 In the case when $a\in S\setminus S^2$ we have
$$X_a= \{b\in S\mid (a,b)\in h \text{ and } b\in S\setminus S^2\}.
$$
Clearly, $(a,b)\in h\iff (a\tau, b\tau)\in h$ and $b\in S\setminus
S^2\iff b\tau\in S\setminus S^2$ for all $a,b\in S$. It follows
that
$$X_{a\tau}=\{b\tau\mid b\in X_a\},
$$
which implies Equation~\eqref{eq:q}.
The inclusion $a\in S^2$ is equivalent to $a\tau\in S^2$. But then
$\vert X_a\vert=\vert X_{a\tau}\vert=1$, which also implies Equation~\eqref{eq:q}.

Suppose now that \eqref{eq:q} holds for certain $\tau\in
{\mathrm{Aut}}T$. Then one can extend $\tau$ to
${\overline{\tau}}\in {Aut{S}}$ as follows.

Fix a collection of sets $I_a$, $a\in T$, and bijections
$f_a:I_a\to X_a$, $a\in T$, satisfying the following conditions:
\begin{itemize}
\item $|I_a|=|X_a|$;
\item $I_a=I_b$ whenever $|X_a|=|X_b|$;
\item $I_a\cap I_b=\emptyset$ whenever $|X_a|\neq|X_b|$;
\item if $a,b\in T$ and $|X_a|=|X_b|$ then $af_a^{-1}=bf_b^{-1}$.
\end{itemize}
Obviously, such  collections $I_a$, $a\in T$, and $f_a$, $a\in T$,
exist.

Consider $x\in S\setminus T$. Since $T$ is a transversal of $\psi$
there is $a\in T$ such that $x\in X_a$. By the hypothesis we have
$|X_a|=|X_{a\tau}|$. Set ${\overline{\tau}}$ on $X_a$ to be the
map from $X_a$ to $X_{a\tau}$, defined via $x\mapsto
xf_a^{-1}f_{a\tau}$. In this way we define a bijection
${\overline{\tau}}$ of $S$ such that
${\overline{\tau}}|_{T}=\tau$. It will be called an {\em
extension} of $\tau$ to $S$. To complete the proof, we are left to
show that ${\overline{\tau}}$ is a homomorphism. Let $x,y\in S$,
$x\in X_a$, $y\in X_b$. Then
$$
(xy){\overline{\tau}}=(ab)\tau=a\tau b\tau= x{\overline{\tau}}
y{\overline{\tau}},
$$
as required.

\end{proof}

Let $\tau \in H$. Of course, ${\overline{\tau}}$, constructed in
the proof of Proposition~\ref{lem:thues}, depends not only on
$\tau$, but also on the sets $I_a$ and the maps $f_a$, so that
$\tau$ may have several extentions to $S$. Fix any extention
${\overline{\tau}}$ of $\tau$.

\begin{lemma}\label{lem:emb} $\tau\mapsto {\overline{\tau}}$ is
an embedding of $H$
into ${\mathrm{Aut}}S$.
\end{lemma}
\begin{proof}
Proof follows directly from the construction of
${\overline{\tau}}$.
\end{proof}

Denote by ${\overline{H}}$ the image of $H$ under the embedding of
$H$ into ${\mathrm{Aut}}S$ from Lemma~\ref{lem:emb}.

\section{Proof of Theorem~\ref{th:main}}

\begin{proposition}\label{pr:action} ${\overline{H}}$ acts on $G$ by automorphisms via
$\pi^{\tau}={\tau}^{-1}\pi{\tau}$, $\tau\in {\overline{H}}, \pi\in
G$.
\end{proposition}
\begin{proof}
Let $\pi\in G$ and $\tau\in {\overline {H}}$. Show first that
$\pi^{\tau}\in G$. Take any $x\in S$. Let $X_a$ be the block which
contains $x$. We consequently have that $x{\tau^{-1}}\in
X_{a\tau^{-1}}$, $x{\tau^{-1}}\pi\in X_{a\tau^{-1}}$ and
$x{\tau^{-1}}\pi\tau\in X_{a\tau^{-1}\tau}=X_a$. Hence,
$x\pi^{\tau}\in X_a$. It follows that $\pi^{\tau}\in G$.

That $\pi\mapsto \pi^{\tau}$ is one-to-one, onto and homomorphic
immediately follows from the definition of this map. We are left
to show that the
 map, which sends $\tau\in{\overline {H}}$ to
$\pi\mapsto \pi^{\tau}\in {\mathrm{Aut}}G$, is homomorphic.
Indeed, $\pi^{\tau_1\tau_2}=(\tau_1\tau_2)^{-1}\pi(\tau_1\tau_2)=
(\pi^{\tau_1})^{\tau_2}$. This completes the proof.
\end{proof}
In the following two lemmas we show that $G$ and ${\overline{H}}$
intersect by the identity automorphism and generate
${\mathrm{Aut}} S$.

\begin{lemma}\label{lem:id_int}
$G\cap{\overline{H}}=id$, where $id$ is the identity automorphism
of $S$.
\end{lemma}

\begin{proof}
The proof follows from the observation that the decomposition
$S=\cup_{a\in T}X_a$ is fixed by each element of $G$, while only
by the identity element of ${\overline{H}}$.
\end{proof}

\begin{lemma}\label{lem:gen}
${\mathrm{Aut}}S={\overline{H}}\cdot G$.
\end{lemma}

\begin{proof} Let $\varphi\in {\mathrm{Aut}}S$. It follows from the definition
of $\psi$ that $\varphi$ maps each $\psi$-class onto some other
$\psi$-class. Define a bijection $\tau: T\to T$ via $a\tau =b$, if
$X_a\varphi=X_b$ and show that $\tau$ is an extendable
automorphism of $T$. It follows immediately from the construction
of $\tau$ that Equation~(\ref{eq:q}) holds, so that  $\tau$ is
extendable by Proposition~\ref{lem:thues}. Let
${\overline{\tau}}\in {\mathrm{Aut}}S$ be an extension of $\tau$.
The construction implies that $\varphi({\overline{\tau}})^{-1}\in
G$.
\end{proof}

Now the proof of Theorem~\ref{th:main} follows from Proposition~\ref{pr:action} and
Lemmas~\ref{lem:id_int} and~\ref{lem:gen}.

%}

\end{document}